\documentclass[12pt, a4paper]{amsart} 

\usepackage{amsmath} \usepackage{amssymb} \usepackage{amsthm}
\usepackage{amscd} 
\usepackage{enumerate}

\def\a{{\mathfrak{a}}} 
\def\b{{\mathfrak{b}}}
\def\F{{\mathbb{F}}} 
 
\def\m{{\mathfrak{m}}}

\def\N{{\mathbb{N}}}

\def\Q{{\mathbb{Q}}}

\def\Hom{{\mathrm{Hom}}}
\def\End{{\mathrm{End}}}

\def\Ann{{\mathrm{Ann}}}
\def\Spec{{\mathrm{Spec\; }}}

\def\p{{\mathfrak p}}
\def\Ass{\operatorname{Ass}}
\def\Max{\operatorname{Max}}
\def\Gor{\operatorname{Gor}}
\def\Ext{\operatorname{Ext}}
\def\Tor{\operatorname{Tor}}
\def\RHom{\operatorname{{\bf R}Hom}}
\def\D{\operatorname{{\mathcal D}}}
\def\Mod{\operatorname{Mod}}

\theoremstyle{plain}
\newtheorem{thm}{Theorem}[section] 
\newtheorem{cor}[thm]{Corollary}
\newtheorem{prop}[thm]{Proposition}

\newtheorem{lem}[thm]{Lemma}
\theoremstyle{definition} 
\newtheorem{defn}[thm]{Definition}
\newtheorem{eg}[thm]{Example} 
\theoremstyle{remark}
\newtheorem{rem}[thm]{Remark}
\newtheorem{ques}[thm]{Question}

\newtheorem*{cl}{Claim}

\newtheorem*{acknowledgement}{Acknowledgments}

\title{$D$-modules over rings with finite\\ F-representation type}
\author{Shunsuke Takagi}
\address{Department of Mathematics, Kyushu University, 
6-10-1 Hakozaki, Higashi-ku, Fukuoka, 812-8581 Japan}
\email{stakagi@math.kyushu-u.ac.jp}
\author{Ryo Takahashi}
\address{Department of Mathematical Sciences, Faculty of Science, Shinshu University, 3-1-1 Asahi, Matsumoto, Nagano 390-8621, Japan}
\email{takahasi@math.shinshu-u.ac.jp}

\begin{document}
\begin{abstract}
Smith and Van den Bergh \cite{SvdB} introduced the notion of finite F-representation type as a characteristic $p$ analogue of the notion of finite representation type. 
In this paper, we prove two finiteness properties of rings with finite F-representation type. 
The first property states that if $R=\bigoplus_{n \ge 0}R_n$ is a Noetherian graded ring with finite (graded) F-representation type, then for every non-zerodivisor $x \in R$, $R_x$ is generated by $1/x$ as a $D_{R}$-module.
The second one states that if $R$ is a Gorenstein ring with finite F-representation type,  then $H_I^n(R)$ has only finitely many associated primes for any ideal $I$ of $R$ and any integer $n$.
We also include a result on the discreteness of F-jumping exponents of ideals of rings with finite (graded) F-representation type as an appendix. 
\end{abstract}

\maketitle
\markboth{S. TAKAGI and R. TAKAHASHI}{$D$-modules over rings with finite F-representation type}

\section*{Introduction}
As a characteristic $p$ analogue of the notion of finite representation type, Smith and Van den Bergh \cite{SvdB} introduced the notion of finite F-representation type (this notion is valid only in prime characteristic).  They then showed that rings with finite F-representation type are suitable for developing the theory of differential operators on the rings. 
Also, rings with finite F-representation type satisfy several finiteness properties. 
For example, Seibert \cite{Se} proved that the Hilbert-Kunz multiplicities of such rings are rational numbers and Yao \cite{Yao} proved that tight closure commutes with localization in such rings. 
Let $k$ be a field,  $R$ be a $k$-algebra and $D_{R/k}$ be the ring of $k$-linear differential operators on $R$. 
In this paper, inheriting the spirit of Smith and Van den Bergh, we pursue finiteness properties of some $D_{R/k}$-modules when $R$ has finite F-representation type. 
We hope that our results will shed a new light on the theory of $D_{R/k}$-modules over singular $k$-algebras $R$ of prime characteristic. 

For every non-zerodivisor $x \in R$, the localization $R_x$ carries a $D_{R/k}$-module structure.  
Hence we can consider an ascending chain of $D_{R/k}$-submodules of $R_x$
$$D_{R/k} \cdot 1/x \subseteq D_{R/k} \cdot 1/x^2 \subseteq D_{R/k} \cdot 1/x^3 \subseteq \cdots \subseteq R_x.$$
Assume that $R=k[x_1, \dots, x_d]$ is a polynomial ring over $k$, and then it is known that this ascending chain stabilizes, that is, $R_x$ is generated by $1/x^i$ for some $i$ as a $D_{R/k}$-module (it was proved by Bernstein \cite{Be} in characteristic zero and by B\o gvad \cite{Bo} in positive characteristic). 
If $k$ is a field of characteristic zero and $-i$ is the smallest integer root of the Bernstein-Sato polynomial $b_x(s)$ of $x$, then $R_x$ is generated by $1/x^i$ as a $D_{R/k}$-module; Walther \cite{Wa} proved that it is not generated by $1/x^j$ for $j<i$.
If $k$ is a field of positive characteristic, then Alvarez-Montaner, Blickle and Lyubeznik \cite{ABL} proved that $R_x$ is generated by $1/x$ as a $D_{R/k}$-module. 
However, when $R$ is not a regular ring, almost nothing is known about the stabilization of the above ascending chain. 
The first ingredient  of this paper is to generalize the result due to Alvarez-Montaner, Blickle and Lyubeznik to the case of rings with finite F-representation type: 
let $R=\bigoplus_{n \ge 0}R_n$ be a Noetherian graded ring with $k=R_0$ a field of positive characteristic. 
If $R$ has finite (graded) F-representation type, then $R_x$ is generated by $1/x$ as a $D_{R/k}$-module (Corollary \ref{generator}).  

Another example of $D_{R/k}$-modules is a local cohomology module $H^n_I(R)$ of $R$.  
In 1990, Huneke \cite{H} raised a problem of whether local cohomology modules of Noetherian rings have finitely many associated primes or not.
It is now known that this problem has a negative answer; several counterexamples were constructed by Singh \cite{Si}, Katzman \cite{K}, and Singh and Swanson \cite{SS}.
However, a lot of affirmative results have been obtained under some additional assumptions.
Huneke and Sharp \cite{HS} and Lyubeznik \cite{L1,L3,L4} proved that if $R$ is a regular ring which either contains a field or is unramified then $H_I^n(R)$ has only finitely many associated primes for any ideal $I$ of $R$ and any integer $n$.
Marley \cite{M} proved that the same assertion holds if $R$ is any Noetherian local ring of Krull dimension at most three.
Thus, determining when the set of associated primes is finite is now an important problem in commutative algebra.
The second ingredient of this paper is to show that if $R$ is a Gorenstein ring with finite F-representation type,  then $H_I^n(R)$ has only finitely many associated primes for any ideal $I$ of $R$ and any integer $n$ (Corollary \ref{assgor}).
This is a generalization of a result of Huneke and Sharp, asserting that the sets of associated primes of local cohomology modules of regular local rings of positive characteristic are finite.
We shall actually prove some more general results (see Proposition \ref{glocus} and Theorem \ref{assthm2}).

In the appendix, we study the discreteness of F-jumping exponents, jumping exponents for generalized test ideals of \cite{HY}. 
These invariants were first studied under the name of F-thresholds in \cite{MTW} and it was shown that they satisfy many of the formal properties of jumping numbers for multiplier ideals. However, the discreteness and rationality of F-jumping exponents were left open in \cite{MTW}. 
Blickle, Musta\c t\v a and Smith \cite{BMS} proved the discreteness and rationality of F-jumping exponents in the case where the ring is an F-finite regular ring of essentially of finite type over a field. 
We prove the discreteness of F-jumping exponents of homogeneous ideals of strongly F-regular graded rings with finite (graded) F-representation type (Theorem \ref{discrete}). 

\section{Rings with finite F-representation type}
Throughout this paper, all rings are Noetherian commutative rings with identity and have prime characteristic $p$. 
For such a ring $R$, the Frobenius map $F:R \to R$ is defined by sending $r$ to $r^p$ for all $r \in R$. 
We often consider not just the Frobenius map but also its iterates $F^e:R \to R$ sending $r$ to $r^{p^e}$. For any $R$-module $M$,  we denote by ${}^e\!M$ the module $M$ with its $R$-module structure pulled back via the $e$-times iterated Frobenius map $F^e$. 
That is, ${}^e\!M$ is just $M$ as an abelian group, but its $R$-module structure is determined by $r \cdot m:=r^{p^e}m$ for all $r \in R$ and $m \in M$. Note that ${}^0\!M=M$. 
If $R=\bigoplus_{n \geq 0}R_n$ is a graded ring and $M$ is a graded $R$-module, then ${}^e\!M$ carries a $\Q$-graded $R$-module structure: we grade ${}^e\!M$ by putting $[{}^e\!M]_\alpha=[M]_{p^e \alpha}$ if $\alpha$ is a multiple of $1/p^e$, otherwise $[M]_{\alpha}=0$.   
Let $I$ be an ideal of $R$. For each $q=p^e$, we denote by $I^{[q]}$ the ideal of $R$ generated by the $q^{\rm th}$ powers of elements of $I$. It is easy to see that $R/I \otimes_R {}^e\!M \cong {}^e\!M/(I \cdot {}^e\!M) \cong {}^e\!(M/I^{[q]}M)$. 
We say that a ring $R$ is \textit{F-pure} if ${}^1\!R$ is a pure extension of $R$ and that $R$ is {\it F-finite} if  ${}^1\!R$ is a finitely generated $R$-module. 

Rings with finite F-representation type were first introduced by Smith and Van den Bergh \cite{SvdB},  under the assumption that the Krull-Schmidt theorem holds for them. 
Yao \cite{Yao} studied these rings in a more general setting.  

\begin{defn}[\textup{\cite[Definition 3.1.1]{SvdB}, \cite[Definition 1.1]{Yao}}]
\ 
\renewcommand{\labelenumi}{(\roman{enumi})}
\begin{enumerate}
\item
Let $R$ be a Noetherian ring of prime characteristic $p$. 
We say that $R$ has \textit{finite F-representation type} by finitely generated $R$-modules $M_1, \dots, M_s$ if for every $e \in \N$,  the $R$-module ${}^e\!R$ is isomorphic to a finite direct sum of the $R$-modules $M_1, \dots, M_s$, that is, there exist nonnegative integers $n_{e1}, \dots, n_{es}$ such that 
$${}^e\!R \cong \bigoplus_{i=1}^s M_i^{\oplus n_{ei}}.$$
We simply say that $R$ has finite F-representation type if there exist finitely generated $R$-modules $M_1, \dots, M_s$ by which $R$ has finite F-representation type. 

\item
Let $R=\bigoplus_{n \geq 0}R_n$ be a Noetherian graded ring of prime characteristic $p$. 
We say that $R$ has \textit{finite graded F-representation type} by finitely generated $\Q$-graded $R$-modules $M_1, \dots, M_s$ if for every $e \in \N$, the $\Q$-graded $R$-module ${}^e\!R$ is isomorphic to a finite direct sum of the $\Q$-graded $R$-modules $M_1, \dots, M_s$, that is,  there exist nonnegative integers $n_{ei}$ and rational numbers $\alpha_{ij}^{(e)}$ for all $1 \le i \le s$ and $1 \le j \le n_{ei}$ such that there exists a $\Q$-graded isomorphism of the form 
$${}^e\!R \cong \bigoplus_{i=1}^s \bigoplus_{j=1}^{n_{ei}}M_i(\alpha_{ij}^{(e)}).$$
We simply say that $R$ has finite graded F-representation type if there exist finitely generated $\Q$-graded $R$-modules $M_1, \dots, M_s$ by which $R$ has finite graded F-representation type. 

\end{enumerate}
\end{defn}

\begin{rem}
(1) When $R$ is reduced, $R$ has finite F-representation type by finitely generated $R$-modules $M_1, \dots, M_s$ if and only if for every $q=p^e$,  the $R$-module $R^{1/q}$ is isomorphic to a finite direct sum of the $R$-modules $M_1, \dots, M_s$.

(2) If the ring $R$ has finite F-representation type, then $R$ is F-finite. 
\end{rem}

We exhibit examples of rings with finite F-representation type. 
\begin{eg}\label{FFRTeg}
\renewcommand{\labelenumi}{(\roman{enumi})}
\begin{enumerate}
\item (\cite[Observation 3.1.2]{SvdB}) 
Let $R$ be an F-finite regular local ring of characteristic $p>0$ (resp.  a polynomial ring $k[X_1, \dots, X_d]$ over a field $k$ of characteristic $p>0$ such that $[k:k^p]<\infty$).
Then $R$ has finite F-representation type (resp. finite graded F-representation type) by the $R$-module $R$. 

\item (\cite[Observation 3.1.3]{SvdB}) Let $R$ be a Cohen-Macaulay local ring (resp. a Cohen-Macaulay graded ring) of prime characteristic $p$ with finite representation type (resp. finite graded representation type). Then $R$ has finite F-representation type (resp. finite graded F-representation type).  

\item Let $(R, \m , k)$ be an Artinian local ring of prime characteristic $p$  (resp. $R=\bigoplus_{n \ge 0} R_n$ be an Artinian graded ring with $k:=R_0$ a field of characteristic $p>0$) such that $[k:k^p] <\infty$. Then $R$ has finite F-representation type (resp. finite graded F-representation type). 
To check this, we may assume that $k$ is a perfect field for simplicity. Since $(R,\m)$ is Artinian, there exists some $q_0=p^{e_0}$ such that $\m^{[q]}=0$ for all $q=p^e \ge q_0$.
Then ${}^e\!R$ is a $k$-vector space for all $e \ge e_0$. 
Since $R$ is F-finite, ${}^e\!R$ is a finitely generated $R$-module.
Hence ${}^e\!R$ is a finite-dimensional $k$-vector space for all $e \ge e_0$. The graded case also follows from a similar argument. 

\item Let $R=k[t^{n_1}, \dots, t^{n_r}]$ be a monomial curve over a field $k$ of characteristic $p>0$ such that $[k:k^p] < \infty$ with $n_i \in \N$. Then $R$ has finite graded F-representation type. 
To check this, we may assume that $k$ is a perfect field and $\gcd(n_1, \dots, n_r)=1$. Let $c$ be the conductor of the semigroup $\langle n_1, \dots, n_r \rangle$, and fix any $q=p^e \ge c$. 
Then $R^{1/q}= \bigoplus_{i=0}^{q-1}(\sum_{v \in S_i}Rt^{v/q})$, where $S_i:=\{v \in \langle n_1, \dots, n_r \rangle \mid v \equiv i \mod q\}$ for $0 \le i \le q-1$. 
Let  $\varphi_i: \sum_{j=0}^{c-1}Rt^{j/q} \to \sum_{v \in S_i}Rt^{v/q}$ be an $R$-module homomorphism defined by $\varphi_i(1)= t^{r_i/q}$, where $r_i$ is the least integer in $S_i$.
Then $\varphi_i$ is an isomorphism for all $i=0, \dots, q-1$. Thus,
$$R^{1/q}= \bigoplus_{i=0}^{q-1}(\sum_{v \in S_i}Rt^{v/q}) \cong (\sum_{j=0}^{c-1}Rt^{j/q})^{\oplus q}.$$

\item (cf.~\cite[Example 2.3.6]{Ka}) Let  $S=k[X_1, \dots, X_d]$ be a polynomial ring over a field $k$ of characteristic $p>0$ such that $[k:k^p]<\infty$ and $I$ be a monomial ideal of $S$. 
Then the quotient $S/I$ has finite graded F-representation type. 

\item (Compare this with \cite[Proposition 3.1.4]{SvdB}) Let $R \hookrightarrow S$ be a finite local homomorphism of Noetherian local rings of prime characteristic $p$ such that $R$ is an $R$-module direct summand of $S$. 
If $S$ has finite F-representation type, so does $R$. 
Indeed, if $S$ has finite F-representation type by finitely generated $S$-modules $N_1, \cdots, N_t$, then for each $e \in \N$, we have a decomposition ${}^e\! S \cong \bigoplus_{j=1}^t N_j^{\oplus m_{ej}}$ with $m_{ej} \in \N$. Put $N=N_1 \oplus \cdots \oplus N_t$.  Since ${}^e\!R$ is a direct summand of ${}^e\!S$, ${}^e\!R$ is a direct summand of a finite direct sum of the finitely generated $R$-module $N$ for all $e \in \N$. 
Then by \cite[Theorem 1.1]{W}, there exist, up to isomorphism, only finitely many indecomposable direct  summands of the $R$-modules ${}^e\!R$'s. Thus, $R$ has finite F-representation type. 

\item (\cite[Proposition 3.1.6]{SvdB}) Let $R=\bigoplus_{n \geq 0}R_n \subseteq S=\bigoplus_{n \geq 0}S_n$ be a Noetherian graded rings with $R_0, S_0$ fields of characteristic $p>0$ such that $R$ is an $R$-module direct summand of $S$.  
Assume in addition that $[S_0:R_0] <\infty$. 
If $S$ has finite graded F-representation type, so does $R$. 
In particular, normal semigroup rings and rings of invariants of linearly reductive groups have finite graded F-representation type. 

\item
Let $R=k[[X,Y,Z]]/(X^2+Y^3+Z^5+XYZ)$ be the simple singularity of type $E_8^4$, where $k$ is an algebraically closed field of characteristic two. 
Then it is well-known that $R$ has finite representation type (see \cite{GK} and its reference for details) and, in particular, $R$ has finite F-representation type by (ii).
On the other hand, we can easily check that $R$ is F-pure but not strongly F-regular (see the paragraph preceding Proposition \ref{test ideal} for the definition of strongly F-regular rings). 

\item 
Let $R$ be the affine cone of a smooth projective curve $X$ of genus $g$ over an algebraically closed field $k$ of characteristic $p>0$.
If $X$ is an elliptic curve, then by \cite{At} and \cite{Ta}, $R$ does not have finite graded F-representation type. 
If $g \ge 2$, then by \cite[Proposition 1.2]{LP}, $R$ does not have finite graded F-representation type. 
\end{enumerate}
\end{eg}

We conclude this section with a couple of questions. 
\begin{ques}
\begin{enumerate}
\item Does every F-finite determinantal ring have finite (graded) F-representation type? 
\item (Brenner) Let $k$ be an algebraically closed field of characteristic $p$. Then does the ring $k[X,Y,Z]/(X^2+Y^3+Z^7)$ have finite F-representation type? 
\end{enumerate}
\end{ques}

\section{Generators of $R_x$ as a $D_R$-module}
For any non-zerodivisor $x \in R$, $R_x$ carries a $D_R$-module structure, where $D_R$ is the ring of differential operators on $R$. In this section, we will study generators of $R_x$ as a $D_R$-module. 
First,  for a given ideal $\a \subseteq R$ and a given module $M$,  we define a descending chain of submodules $\{I_e(\a, M)\}_{e \ge 0}$ of $M$.  
\begin{defn}\label{idealdef}
Let $R$ be a Noetherian ring of prime characteristic $p$ and $M$ be an $R$-module. 
For any integer $e \ge 0$ and any ideal $\a$ of $R$, the $R$-submodule $I_e(\a, M) \subseteq M$ is defined to be
$$I_e(\a, M):=\Hom_R({}^e\!R, M)\cdot \a=\{\phi(x) \mid \phi \in \Hom_R({}^e\!R,M), x \in \a\} \subseteq M.$$
When $R$ is reduced, we can think of $I_e(\a, M)$ as 
$$I_e(\a, M)=\Hom_R(R^{1/p^e}, M)\cdot \a^{1/p^e} \subseteq M.$$
\end{defn}

\begin{rem}
When $R$ is an F-finite regular ring, the ideal $I_e(x, R)$ coincides with the ideal $I_e(x)$ in \cite{ABL} and the ideal $I_e(\a^r, R)$ does with the ideal $(\a^r)^{[1/p^e]}$ in \cite{BMS}.  
The reader is also referred to Proposition \ref{test ideal}.  
\end{rem}

We list up basic properties of  $I_e(\a, M)$ which we will need later. 
Before doing so, we recall the definition of F-pure pairs. The pair $(R,x)$ of an F-finite Noetherian reduced ring of characteristic $p>0$ and a non-zerodivisor $x \in R$ is said to be \textit{F-pure} if the inclusion $x^{(q-1)/q}R \hookrightarrow R^{1/q}$ splits as an $R$-linear map for all $q=p^e$.  For example, when $R$ is Gorenstein, $(R,x)$ is an F-pure pair if and only if $R/(x)$ is an F-pure ring (see \cite[Remark 4.10]{HW}).

\begin{lem}\label{idealbasic}
Let $R,M,e,\a$ be the same as in Definition \ref{idealdef}. In addition, suppose that $R$ is F-finite. 
\renewcommand{\labelenumi}{$(\arabic{enumi})$}
\begin{enumerate}
\item For ideals $\a \subseteq \b$ in $R$, $I_e(\a, M) \subseteq I_e(\b,M)$. 
\item $I_e(\a, M) \supseteq I_{e+1}(\a^{[p]}, M)$. 
If $R$ is F-pure, then $I_e(\a, M) = I_{e+1}(\a^{[p]}, M)$. 
\item For every $x \in R$ and every $q=p^e$, $I_{e+1}(x^{pq-1}, M) \subseteq I_e(x^{q-1}, M).$ 
Moreover, if the pair $(R,x)$ is F-pure, then $I_e(x^{p^e-1}, M)=M$ for every integer $e \ge 0$. 
\end{enumerate}
\end{lem}
\begin{proof}
(1) obvious. 

(2) Each $R$-linear map $\phi_{e+1}:{}^{e+1}\!R \to M$ induces an $R$-linear map $\phi_{e}:{}^e\!R \to M$ such that ${\phi}_{e}(x)=\phi_{e+1}(x^p)$ for all $x \in R$. 
Thus, $\phi_{e+1}(\a^{[p]})=\phi_e(\a) \subseteq I_{e}(\a, M)$. 
If $R$ is F-pure, then there exists an ${}^e\!R$-linear map $f_e:{}^{e+1}\!R \to {}^e\!R$ sending $1$ to $1$. 
For any $\psi_e \in \Hom_R({}^{e}\!R,M)$,  composing with $f_e$, we have an $R$-linear map $\psi_{e+1}:{}^{e+1}\!R \to M$ such that $\psi_{e+1}(x^p)=\psi_e(x)$ for all $x \in R$. 
Therefore, $\psi_e(\a)=\psi_{e+1}(\a^{[p]}) \subseteq  I_{e+1}(\a^{[p]}, M)$.

(3) Since $x^{pq-1}=(x^{q-1})^px^{p-1}$, by (1) and (2), one has 
$$I_{e+1}(x^{pq-1}, M) \subseteq I_{e+1}((x^{q-1})^p, M) \subseteq I_e(x^{q-1}, M).$$
If $(R,x)$ is F-pure, then for every $e \in \N$, there exists an $R$-linear map $f:{}^{e}\!R \to R$ sending $x^{p^e-1}$ to $1$. 
For each $m \in M$, $f$ induces an $R$-linear map ${}^e\!R \to M$ sending $x^{p^e-1}$ to $m$. This implies that $I_e(x^{p^e-1}, M)=M$. 
\end{proof}

Now we recall the definition of rings of differential operators. The reader is referred to \cite{GD} and \cite{Ye} for details.
\begin{defn}\label{Ddef}
Let $k$ be any commutative ring and $R$ be any commutative $k$-algebra. 
The ring of $k$-linear differential operators $D_{R/k}$ is a subring of the $k$-linear endomorphism ring $\mathrm{End}_k(R)$ of $R$. The ring $D_{R/k}$ is the union $\bigcup_{n \ge 0} D^n_{R/k}$, where $D_{R/k}$ are defined inductively as follows:
\begin{align*}
D^{-1}_{R/k}&=0\\
D^{n+1}_{R/k}&=\{\delta \in \Hom_k(R,R) \mid [\delta, r] \in D^n_{R/k}\textup{ for all $r \in R$}\}.
\end{align*}
Here $[\delta, r]=\delta \circ r -r \circ \delta$ denotes the commutator. 
\end{defn}

\begin{prop}[\textup{\cite[1.4.9]{Ye}}]\label{F-finite}
Let the notation be the same as in Definition \ref{Ddef} and assume that $R$ is a ring of prime characteristic $p$. 
Then $D_{R/k}$ is a subring of $\bigcup_{q=p^e}\mathrm{End}_{R^{q}}(R)$.
If $R$ is F-finite and $k$ is a perfect field of characteristic $p>0$, then
$$D_{R/k}=\bigcup_{q=p^e}\mathrm{End}_{R^{q}}(R).$$
\end{prop}
The above proposition implies that $D_{R/k}$ is the same for every perfect field $k \subseteq R$.  Consequently, we write simply $D_R$ instead of $D_{R/k}$ when $k$ is a perfect field. 
We set $D_R^{(e)}=\End_{R^q}(R)$.

\begin{eg}[\textup{\cite[Ex. 5.1]{L2}}]
Let  the notation be the same as in Definition \ref{Ddef} and assume that $R$ is a ring of  prime characteristic $p$. 
\renewcommand{\labelenumi}{(\roman{enumi})}
\begin{enumerate}
\item
$R$ itself is a $D_{R/k}$-module with its natural $D_{R/k}$-action. 
\item
A localization $R_S$ with respect to any multiplicatively closed subset $S \subset R$ carries a unique $D_{R/k}$-module structure such that the natural localization map  $R \to R_S$ is $D_{R/k}$-linear:  the operator $\delta \in D_{R/k}$ is $R^{q}$-linear for $q=p^e \gg 0$ by Proposition \ref{F-finite}. For such $q$, we define 
$$\delta\left(\frac{x}{s}\right)=\delta\left(\frac{s^{q-1}x}{s^{q}}\right)=\frac{\delta(s^{q-1}x)}{s^{q}}.$$
The well-definedness and uniqueness are easily checked. 
\item
Let $I$ be an ideal of $R$ and $n$ be a nonnegative integer. 
Then by virtue of (ii), the local cohomology module $H^i_I(R)$ has a $D_{R/k}$-module structure, because $H^i_I(R)$ is the $i^{\rm th}$ cohomology module of a \v Cech complex. 
\end{enumerate}
\end{eg}

\begin{prop}[\textup{cf. \cite[Proposition 3.5]{ABL}}]\label{equiv}
Let $x$ be a non-zerodivisor of a Noetherian ring $R$ of prime characteristic $p$. 
Suppose that $R$ has finite F-representation type by finitely generated $R$-modules $M_1, \dots, M_s$. If the descending chain of submodules in $M_i$
$$I_1(x^{p-1}, M_i) \supseteq I_2(x^{p^2-1}, M_i) \supseteq I_3(x^{p^3-1}, M_i) \supseteq \dots$$
stabilizes at $e$ $($that is,  $I_e(x^{p^e-1}, M_i)=I_{e+1}(x^{p^{e+1}-1}, M_i)=\dots)$ for all $i=1, \dots, s$, then there exists $\delta \in D^{(e+1)}_R$ such that $\delta(1/x)=1/x^p$. 
When $R$ is F-pure, the converse also holds true. 
\end{prop}
\begin{proof}
First, we assume that $I_e(x^{p^e-1}, M_i)=I_{e+1}(x^{p^{e+1}-1}, M_i)$ for all $i=1, \dots, s$. 
Then by Lemma \ref{idealbasic} (2), $I_{e+1}(x^{p^{e+1}-p}, M_i)=I_{e+1}(x^{p^{e+1}-1}, M_i)$. 
Since $R$ has finite F-representation type by $M_1, \dots, M_s$, 
we have an decomposition 
$$\End_R({}^{e+1}\!R) \cong \bigoplus_{i=1}^s \Hom_R({}^{e+1}\!R, M_i)^{\oplus n_{e+1i}}.$$
Thus, $\End_R({}^{e+1}\!R) \cdot x^{p^{e+1}-p}=\End_R({}^{e+1}\!R) \cdot x^{p^{e+1}-1}$. 
This implies that $D^{(e+1)}_R \cdot  x^{p^{e+1}-p}=D^{(e+1)}_R \cdot  x^{p^{e+1}-1}$. In particular, there exists $\delta \in D^{(e+1)}_R$ such that $\delta(x^{p^{e+1}-1})=x^{p^{e+1}-p}$. 
Dividing this equality by $x^{p^{e+1}}$ gives that $\delta(1/x)=1/x^p$. 

Next, we will prove the converse. Just reversing the above argument, we have that $I_{e+1}(x^{p^{e+1}-p}, M_i)=I_{e+1}(x^{p^{e+1}-1}, M_i)$ for all $i=1, \dots, s$.  
Since $R$ is F-pure, by Lemma \ref{idealbasic} (2) again, $I_{e}(x^{p^{e}-1}, M_i)=I_{e+1}(x^{p^{e+1}-1}, M_i)$ for all $i=1, \dots, s$.
\end{proof}

\begin{ques}
Let $R$ be a Noetherian ring with finite F-representation type by finitely generated $R$-modules $M_1, \dots, M_s$ and let $x \in R$ be an element. Does the descending chain of submodules in $M_i$
$$I_1(x^{p-1}, M_i) \supseteq I_2(x^{p^2-1}, M_i) \supseteq I_3(x^{p^3-1}, M_i) \supseteq \dots$$
stabilize for every $i=1, \dots, s$?
\end{ques}

In general, the above question is open. However, when $R$ is a Noetherian graded ring, we can give an affirmative answer to this question. 
\begin{thm}\label{stabilize}
Let $R=\bigoplus_{n \geq 0}R_n$ be a Noetherian graded ring with $R_0$ a field of characteristic $p>0$. 
Assume that $R$ has finite graded F-representation type by finitely generated $\Q$-graded $R$-modules $M_1, \dots, M_s$. 
Then for every $x \in R$ and for every $i=1, \dots, s$, the descending chain of submodules in $M_i$
$$I_1(x^{p-1}, M_i) \supseteq I_2(x^{p^2-1}, M_i) \supseteq I_3(x^{p^3-1}, M_i) \supseteq \dots$$
stabilizes. 
\end{thm}
\begin{proof} 
By assumption, for each $e \in \N$, we have a degree-preserving isomorphism ${}^e\!R \cong \bigoplus_{i=1}^s \bigoplus_{j=1}^{n_{ei}}M_i(\alpha_{ij}^{(e)})$. 
Without loss of generality, we may assume that the grading on the $M_i$ is normalized so that $[M_i]_0 \ne 0$ and $[M_i]_{\alpha}=0$ for any $\alpha<0$. 
Since $[{}^e\!R]_{\alpha}=0$ for any $\alpha<0$, we easily see that $\alpha_{ij}^{(e)} \le 0$ for all $e \in \N$, $i=1, \dots, s$ and $j=1, \dots, n_{ei}$. 
\begin{cl}
The set $\{\alpha_{ij}^{(e)}\}_{e \in \N, 1 \le i \le s,1\le j \le n_{ei}}$ is bounded. 
\end{cl}
\begin{proof}[Proof of Claim]
This claim follows from \cite[Lemma 3.1.5]{SvdB}, but we include the proof for the sake of completeness. 

First note that the Krull-Schmidt theorem holds for $R$, and we may assume that the $M_i$ are indecomposable $\Q$-graded $R$-modules. For each $1 \le k \le s$, by assumption, we have a decomposition of ${}^1\!M_k$ as a $\Q$-graded $R$-module 
$${}^1\!M_k \cong \bigoplus_{i=1}^s \bigoplus_{j=1}^{r_{ki}}M_i(\beta_{ij}^{(k)})$$
with $\beta_{ij}^{(k)} \in \Q_{\le 0}$. 
Set $C_1=\min_{k,i,j}\{\beta_{ij}^{(k)}\}$ and $C_2=\min_{i,j}\{\alpha_{ij}^{(1)}\}$, and then $C_1$ and $C_2$ are non-positive rational numbers independent of $e$. 
Since 
$$\bigoplus_{i=1}^s \bigoplus_{j=1}^{n_{e+1i}}M_i(\alpha_{ij}^{(e+1)}) \cong {}^{e+1}\!R \cong \bigoplus_{k=1}^s \bigoplus_{l=1}^{n_{ek}}{}^1\!M_k(\alpha_{kl}^{(e)}/{p}), $$
one has $\alpha_{ij}^{(e+1)} \ge  C_1+\min_{k,l}\{\alpha_{kl}^{(e)}/p\}$ for every $e \in \N$. 
Thus, we obtain by an easy induction that for all $e \in \N$, $i=1, \dots, s$ and $j=1, \dots, n_{ei}$, 
$$\alpha_{ij}^{(e)} \ge C_1+\frac{C_1}{p}+\dots+\frac{C_1}{p^{e-1}}+\frac{C_2}{p^e} \ge \frac{pC_1}{p-1}+C_2.$$
\end{proof}

Let $x \in R$. We can write
$x^{p^e-1}=\sum_{i=1}^s\sum_{j=1}^{n_{ei}}m_{ij}^{(e)}$ with $m_{ij}^{(e)} \in M_i(\alpha_{ij}^{(e)})$ via the identification of ${}^e\!R$ with $\bigoplus_{i=1}^s \bigoplus_{j=1}^{n_{ei}}M_i(\alpha_{ij}^{(e)})$.  
Then one has 
$$\frac{p^e-1}{p^e} \deg_R(x)=\deg_{{}^e\!R}(x^{p^e-1}) \ge \deg_{M_i(\alpha_{ij}^{(e)})}(m_{ij}^{(e)})$$ for all $e \in \N$, $i=1, \dots, s$ and $j=1, \dots, n_{ei}$, and in particular, $\deg_R(x)>\deg_{M_i(\alpha_{ij}^{(e)})}(m_{ij}^{(e)})$.
On the other hand, by the above claim, we can take a nonnegative rational number $C$ such that $-\alpha_{ij}^{(e)} \le C$ for all $e \in \N$, $i=1, \dots, s$ and $j=1, \dots, n_{ei}$. Thus, for each $1 \le k \le s$, 
\begin{align*}
I_e(x^{p^e-1}, M_k)
&=\sum_{i=1}^s\sum_{j=1}^{n_{ei}}\Hom_{R}(M_i(\alpha_{ij}^{(e)}), M_k) \cdot m_{ij}^{(e)}\\
&\subseteq \sum_{i=1}^s\sum_{j=1}^{n_{ei}}\Hom_{R}(M_i(\alpha_{ij}^{(e)}), M_k) \cdot [M_i(\alpha_{ij}^{(e)})]_{<\deg_R(x)}\\
&\subseteq \sum_{i=1}^s \Hom_{R}(M_i, M_k) \cdot [M_i]_{<\deg_R(x)+C}.
\end{align*}
Let $\phi_{ik}^{(1)}, \dots, \phi_{ik}^{(f_{ik})}$ be a system of homogeneous generators of $\Hom_{R}(M_i, M_k)$ as a $\Q$-graded $R$-module for every $i=1, \dots, s$. 
Then one can choose a rational number $d_k$ such that $\phi_{ik}^{(l)}([M_i]_{<\deg(x)+C}) \subseteq [M_k]_{\le d_k}$ for all $i=1, \dots, s$ and $l=1, \dots, f_{ik}$. 
Since $[M_k]_{\le d_k}$ is a finite-dimensional vector space over $R_0$, the descending chain of subspaces in $[M_k]_{\le d_k}$
$$I_1(x^{p-1}, M_k)\cap [M_k]_{\le d_k} \supseteq I_2(x^{p^2-1}, M_k) \cap [M_k]_{\le d_k} \supseteq I_3(x^{p^3-1}, M_k)\cap [M_k]_{\le d_k} \supseteq \dots$$
stabilizes. This implies that the descending chain of submodules in $M_k$
$$I_1(x^{p-1}, M_k) \supseteq I_2(x^{p^2-1}, M_k) \supseteq I_3(x^{p^3-1}, M_k) \supseteq \dots$$ stabilizes, because $I_e(x^{p^e-1}, M_k)\cap [M_k]_{\le d_k}$ generates $I_e(x^{p^e-1}, M_k)$ as an $R$-module for every $e \in \N$.  
\end{proof}

As a corollary of the above theorem, we can generalize the result due to Alvarez-Montaner, Blickle and Lyubeznik (\cite[Corollary 3.8]{ABL}) to the case of rings with finite F-representation type. 

\begin{cor}\label{generator}
Let $R=\bigoplus_{n \geq 0}R_n$ be a Noetherian graded ring with $k:=R_0$ a field of characteristic $p>0$ and $K$ be the perfect closure of $k$. 
If $R_K:=R \otimes_{k} K$ has finite graded F-representation type, then for any non-zerodivisor $x \in R$, $R_x$ is generated by $1/x$ as a $D_{R/k}$-module. 
\end{cor}
\begin{proof}
We employ the same strategy as the proof of \cite[Corollary 3.8]{ABL}. 
First we prove the following claim. 
\begin{cl}
For any non-zerodivisor $x \in R$, there exists a differential operator $\delta \in D_{R/k}$ such that $\delta(1/x)=1/x^p$. 
\end{cl}
\begin{proof}[Proof of Claim]
By Proposition \ref{equiv} and Theorem \ref{stabilize}, there exists a differential operator $\delta_K \in  D_{R_K/K}$ such that $\delta_K(1/x)=1/x^p$. 
Since $D_{R_K/K}=D_{R/k} \otimes_k K$ by \cite{GD}, we can write $\delta_K=\sum_i c_i \delta_i$ where $c_i \in K$ and $\delta_i \in D_{R/k}$. 
Let $\{\alpha_{\lambda}\}$ be the basis of $R$ over $k$. 
Multiplying both sides of the equation $\sum_i c_i \delta_i(1/x)=1/x^p$ by a sufficiently large power of $x$,  we have an equation $\sum_j \sum_i c_i d_{i\lambda_j}\alpha_{\lambda_j}=0$ in $R_K$ with $d_{i\lambda_j} \in k$.  This means that a system of finitely many linear equations $\{\sum_i X_i d_{i\lambda_j}=0\}_j$ with coefficients in $k$ has a solution $X_i=c_i$ for all $i$. 
By elementary linear algebra, it should have a solution in $k$. Thus, by reversing the above argument, there exists a differential operator $\delta \in D_{R/k}$ such that $\delta(1/x)=1/x^p$. 
\end{proof}

Let $x$ be a non-zerodivisor of $R$.  
Then it follows from the repeated applications of the above claim that
$$D_{R/k} \cdot 1/x=D_{R/k} \cdot 1/x^p=D_{R/k} \cdot 1/x^{p^2}= \cdots =D_{R/k} \cdot 1/x^{p^e}$$
for every $e \in \N$. 
Since the set $\{1/x^{p^e}\}_{e \in \N}$ generates $R_x$ as a $D_{R/k}$-module, we obtain the assertion.  
\end{proof}

\begin{rem}
When the pair $(R,x)$ of an F-finite Noetherian reduced ring $R$ of characteristic $p>0$ and a non-zerodivisor $x \in R$ is F-pure, $R_x$ is always generated by $1/x$ as a $D_R$-module (without assuming that $R$ has finite F-representation type). Indeed, by Lemma \ref{idealbasic} (3), one has $D^{(1)} \cdot f^{p-1}=R$. This means that there exists $\delta \in D^{(1)}_R$ such that $\delta(1/x)=1/x^p$.  
\end{rem}

\section{Associated primes of local cohomology modules}
In this section, we will study the number of associated primes of local cohomology modules when the base ring has finite F-representation type. 
We begin with the following well-known lemma. 

\begin{lem}\label{asslem}
Let $R$ be a Noetherian ring.
\begin{enumerate}[\rm (1)]
\item
Let $\{M_\lambda\}_{\lambda\in\Lambda}$ be a direct system of $R$-modules.
Then $\Ass_R(\varinjlim M_\lambda)$ is contained in $\bigcup_{\lambda\in\Lambda}\Ass_RM_\lambda$.
\item
Suppose that $R$ is a ring of prime characteristic.
Let $M$ be an $R$-module.
Then $\Ass_RM=\Ass_R({}^e\!M)$ for every $e\in \N$.
\end{enumerate}
\end{lem}

\begin{proof}
(1) Let $\p$ be a prime ideal of $R$.
Assume that $\p$ is not an associated prime of $M_\lambda$ for any $\lambda\in\Lambda$.
Then we have $\Hom_{R_\p}(\kappa(\p),(M_\lambda)_\p)=0$ for every $\lambda\in\Lambda$, where $\kappa(\p)=R_\p/\p R_\p$.
Hence $\Hom_{R_\p}(\kappa(\p),(\varinjlim M_\lambda)_\p)\cong\varinjlim\Hom_{R_\p}(\kappa(\p),(M_\lambda)_\p)=0$.
Thus $\p$ is not an associated prime of $\varinjlim M_\lambda$.

(2) Let $\p$ be a prime ideal of $R$.
To show that $\p\in\Ass_RM$ if and only if $\p\in\Ass_R({}^e\!M)$, by localizing at $\p$, it is enough to prove that $\m\in\Ass_RM$ if and only if $\m\in\Ass_R({}^e\!M)$ for a module $M$ over a local ring $(R,\m)$ of prime characteristic.
This is equivalent to saying that $\varGamma_\m(M)\ne 0$ if and only if $\varGamma_\m({}^e\!M)\ne 0$, but this is straightforward.
\end{proof}

Now let us prove the following theorem concerning local cohomology modules of a dualizing complex.

\begin{thm}\label{assmain}
Let $R$ be a Noetherian ring with finite F-representation type by finitely generated $R$-modules $M_1,\dots,M_t$.
Let $D$ be a dualizing complex of $R$, $I$ an ideal of $R$ and $n$ an integer.
Then one has
$$
\Ass_R H_I^n(D)\subseteq\bigcup_{i=1}^t\Ass_R\Ext_R^n(M_i/IM_i,D).
$$
In particular, $H_I^n(D)$ has only finitely many associated primes.
\end{thm}

\begin{proof}
We can assume without loss of generality that $R$ is a local ring.
By definition, for each $e\in \N$, we have an isomorphism of $R$-modules of the form
$$
{}^e\!R\cong\bigoplus_{i=1}^tM_i^{\oplus m_{ei}}.
$$
We denote by $\D(R)$ the derived category of $R$ and by $D_S$ a normalized dualizing complex of a ring $S$.
Note that the additive exact functor ${}^e\!(-):\Mod R\to\Mod{}^e\!R$ extends to an additive exact functor from $\D(R)$ to $\D({}^e\!R)$.
Since $R$ is F-finite, ${}^e\!R$ is a module-finite $R$-algebra.
Hence we have an isomorphism $D_{{}^e\!R}\cong\RHom_R({}^e\!R,D_R)$ in $\D({}^e\!R)$ by \cite[Theorem (3.9)]{Sh}.
Consequently we have the following isomorphisms in $\D(R)$:
\begin{align*}
{}^e\!\RHom_R(R/I^{[p^e]},D_R) 
& \cong \RHom_{{}^e\!R}({}^e\!(R/I^{[p^e]}),{}^e\!(D_R)) \\
& \cong \RHom_{{}^e\!R}({}^e\!R/I\,{}^e\!R,\RHom_R({}^e\!R,D_R)) \\
& \cong \RHom_R({}^e\!R/I\,{}^e\!R,D_R) \\
& \cong \bigoplus_{i=1}^t\RHom_R(M_i/IM_i,D_R)^{\oplus m_{ei}}.
\end{align*}
Therefore there exists an isomorphism 
$${}^e\!\Ext_R^n(R/I^{[p^e]},D_R)\cong\bigoplus_{i=1}^t\Ext_R^n(M_i/IM_i,D_R)^{\oplus m_{ei}}$$ of $R$-modules, and we obtain 
$$\Ass_R\Ext_R^n(R/I^{[p^e]},D_R)=\Ass_R{}^e\!\Ext_R^n(R/I^{[p^e]},D_R)=\bigcup_{i=1}^t\Ass_R\Ext_R^n(M_i/IM_i,D_R)$$ 
by Lemma \ref{asslem} (2). 
Finally, applying Lemma \ref{asslem} (1), we get the desired inclusion 
$$\Ass_RH_I^n(D_R)=\Ass_R(\varinjlim_e\Ext_R^n(R/I^{[p^e]},D_R))\subseteq\bigcup_{i=1}^t\Ass_R\Ext_R^n(M_i/IM_i,D_R).$$
\end{proof}

As a direct corollary of Theorem \ref{assmain}, we obtain the following result.

\begin{cor}\label{assgor}
Let $R$ be a Cohen-Macaulay ring with finite F-representation type by finitely generated $R$-modules $M_1,\dots,M_t$, admitting a canonical module $\omega_R$. 
Let $I$ be an ideal of $R$ and $n$ an integer.
Then one has
$$
\Ass_R H_I^n(\omega_R)\subseteq\bigcup_{i=1}^t\Ass_R\Ext_R^n(M_i/IM_i,\omega_R).
$$
In particular, $H_I^n(\omega_R)$ has only a finite number of associated primes.
\end{cor}

\begin{rem}\label{assrem}
(1) We can prove a little more general form of Corollary \ref{assgor}: let $(R,\m,k)$ be a Cohen-Macaulay local ring of prime characteristic $p$ with canonical module $\omega_R$ and $K$ be an extension field of  the residue field $k$. Assume that $\widehat{R} \widehat{\otimes}_k K$ has finite F-representation type, where $\widehat{R}$ is the $\m$-adic completion of $R$ and $\widehat{\otimes}_k$ denotes the complete tensor product over $k$. Let $I$ be an ideal of $R$ and $n$ an integer. Then
$$
\Ass_R H_I^n(\omega_R)\subseteq\bigcup_{i=1}^t\Ass_R\Ext_R^n(M_i/IM_i,\omega_R).
$$

(2)
Singh and Swanson \cite{SS} constructed an example of a strongly F-regular hypersurface singularity $R$ which is a unique factorization domain such that $H_I^n(R)$ has infinitely many associated primes for some ideal $I$ of $R$ and some integer $n$.
Corollary \ref{assgor} especially says that this ring does not have finite F-representation type.
\end{rem}

Since every F-finite regular local ring has finite F-representation type by the $R$-module $R$, Corollary \ref{assgor} (Remark \ref{assrem} (1)) yields the following result due to Huneke and Sharp \cite[Corollary 2.3]{HS}.

\begin{cor}[Huneke-Sharp]
Let $R$ be a regular local ring of positive characteristic, $I$ an ideal and $n$ an integer.
Then
$$
\Ass_R H_I^n(R)\subseteq\Ass_R\Ext_R^n(R/I,R).
$$
In particular, there exist only finitely many associated primes of $H_I^n(R)$.
\end{cor}

\begin{rem}
As an analogy of another result due to Huneke and Sharp \cite[Theorem 2.1]{HS}, one might expect that the Bass numbers of a local cohomology module should be finite if the ring has finite F-representation type. 
This is, however, false in general. For example, set $R=\F_p[x,y,z,w]/(xz-yw)$, $\m=(x,y,z,w)$ and $I=(x,y)$. 
Since $R$ is a normal Gorenstein semigroup ring, by Example \ref{FFRTeg} (vii), $R$ has finite F-representation type. 
Hartshorne \cite[Section 3]{Ha} showed that $H^2_I(R) \ne 0$ and the $0^{\rm th}$ Bass number $\mu_0(\m,H^2_I(R))$ is infinite. 
\end{rem}

Let $R$ be a Gorenstein ring with finite F-representation type.
Then Corollary \ref{assgor} implies that the set $\Ass_RH_I^n(R)$ is finite.
As for the finiteness of the set $\Ass_RH_I^n(R)$, we can indeed prove a little stronger statement in the case where $R$ is a semilocal ring.
We denote by $\Gor R$ the Gorenstein locus of a ring $R$, namely the set of all prime ideals $\p$ of $R$ such that $R_\p$ is a Gorenstein local ring.

\begin{prop}\label{glocus}
Let $R$ be a Noetherian ring with finite F-representation type by finitely generated $R$-modules $M_1,\dots,M_t$.
Let $I$ be an ideal of $R$ and $n$ an integer.
Then
$$
\Ass_RH_I^n(R)\cap\Gor R\subseteq\bigcup_{i=1}^t\Ass_R\Ext_R^n(M_i/IM_i,R).
$$
In particular, $\Ass_RH_I^n(R)$ is a finite set if $R$ is a semilocal ring such that $R_\p$ is Gorenstein for each nonmaximal prime ideal $\p$.
\end{prop}

\begin{proof}
Let $\p\in\Ass_RH_I^n(R)\cap\Gor R$.
Then $R_\p$ is a Gorenstein local ring, and $\p R_\p$ is in $\Ass_{R_\p}H_{IR_\p}^n(R_\p)$, which is contained in the set $\bigcup_{i=1}^t\Ass_{R_\p}\Ext_{R_\p}^n((M_i)_\p/I(M_i)_\p,R_\p)$ by Corollary \ref{assgor}.
Therefore $\p$ is in $\bigcup_{i=1}^t\Ass_R\Ext_R^n(M_i/IM_i,R)$.
\end{proof}

Next, we consider the finiteness of the sets of associated primes of local cohomology modules of Cohen-Macaulay rings with finite F-representation type.
We denote by $(-)^\dag$ the canonical dual $\Hom_R(-,\omega_R)$, where $\omega_R$ is a canonical module of a Cohen-Macaulay ring $R$.

\begin{lem}\label{assthm}
Let $R$ be a Cohen-Macaulay ring with finite F-representation type by finitely generated $R$-modules $M_1,\dots,M_t$, admitting a canonical module $\omega_R$.
Suppose that the local ring $R_\p$ is Gorenstein for every nonmaximal prime ideal $\p$ of $R$.
Let $I$ be an ideal of $R$ and $n$ an integer.
Then
$$
\Ass_RH_I^n(R)\subseteq\bigcup_{i=1}^t\Ass_R\Ext_R^n(M_i^\dag/IM_i^\dag,\omega_R).
$$
\end{lem}

\begin{proof}
We may assume that $R$ is a local ring.
Moreover, we can also assume that $n\le d:=\dim R$, since $H_I^i(R)=0$ for every $i>d$.
Denote by $\omega_S$ a canonical module of a Cohen-Macaulay ring $S$.
Note that there are isomorphisms 
\begin{align*}
\RHom_R(R/I^{[p^e]},R)&\cong\RHom_R(R/I^{[p^e]},\RHom_R(\omega_R,\omega_R))\\
&\cong\RHom_R(R/I^{[p^e]}\otimes_R^{\bf L}\omega_R,\omega_R).
\end{align*}
Hence there is a spectral sequence
$$
E_2^{ij}=\Ext_R^i(\Tor_j^R(R/I^{[p^e]},\omega_R),\omega_R)\ \Rightarrow\ H^{i+j}=\Ext_R^{i+j}(R/I^{[p^e]},R).
$$
Since $R$ is Gorenstein on the punctured spectrum of $R$, we have $(\omega_R)_\p\cong\omega_{R_\p}\cong R_\p$ for every $\p\in\Spec R \setminus \{\m\}$.
Hence the $R$-module $\Tor_j^R(R/I^{[p^e]},\omega_R)$ has finite length for $j>0$.
It is seen that $E_2^{ij}=0$ if either $i<d$ and $j>0$ or $i>d$.
Thus, from the spectral sequence, we obtain an isomorphism
$$
\Ext_R^n(R/I^{[p^e]},R)\cong\Ext_R^n(\omega_R/I^{[p^e]}\omega_R,\omega_R)
$$
of $R$-modules.
For all $e \in \N$, there is an isomorphism of the form ${}^e\!R\cong\bigoplus_{i=1}^tM_i^{\oplus m_{ei}}$.
Applying $(-)^\dag$ to this isomorphism gives 
$${}^e\!(\omega_R)\cong\omega_{{}^e\!R}\cong\Hom_R({}^e\!R,\omega_R)\cong\bigoplus_{i=1}^t(M_i^\dag)^{\oplus m_{ei}}.$$
Similarly to the arguments in the proof of Theorem \ref{assmain}, we obtain isomorphisms
\begin{align*}
{}^e\!\Ext_R^n(R/I^{[p^e]},R) 
& \cong {}^e\!\Ext_R^n(\omega_R/I^{[p^e]}\omega_R,\omega_R) \\
& \cong \Ext_{{}^e\!R}^n({}^e\!(\omega_R)/I\,{}^e\!(\omega_R),{}^e\!(\omega_R)) \\
& \cong \Ext_{{}^e\!R}^n({}^e\!(\omega_R)/I\,{}^e\!(\omega_R),\RHom_R({}^e\!R,\omega_R)) \\
& \cong \Ext_R^n({}^e\!(\omega_R)/I\,{}^e\!(\omega_R),\omega_R) \\
& \cong \bigoplus_{i=1}^t\Ext_R^n(M_i^\dag/IM_i^\dag,\omega_R)^{\oplus m_{ei}},
\end{align*}
and we get the desired inclusion $\Ass_RH_I^n(R)\subseteq\bigcup_{i=1}^t\Ass_R\Ext_R^n(M_i^\dag/IM_i^\dag,\omega_R)$.
\end{proof}

Using Lemma \ref{assthm}, we can show the following result.
Compare it with Proposition \ref{glocus}.

\begin{thm}\label{assthm2}
Let $R$ be a Cohen-Macaulay ring with finite F-representation type by finitely generated $R$-modules $M_1,\dots,M_t$, admitting a canonical module $\omega_R$.
Suppose that the local ring $R_\p$ is Gorenstein for any prime ideal $\p$ of $R$ with $\dim R/\p\ge 2$.
Let $I$ be an ideal of $R$ and $n$ an integer.
Then one has
$$
\Ass_RH_I^n(R)\subseteq\bigcup_{i=1}^t\Ass_R\Ext_R^n(M_i^\dag/IM_i^\dag,\omega_R)\cap\Max R.
$$
In particular, $\Ass_RH_I^n(R)$ is a finite set if $R$ is semilocal.
\end{thm}

\begin{proof}
Let $\p$ be a nonmaximal associated prime of $H_I^n(R)$.
From the assumption, we easily see that $R_\p$ is Gorenstein on the punctured spectrum.
Lemma \ref{assthm} implies that the set $\Ass_{R_\p}H_{I_\p}^n(R_\p)$ is contained in $\bigcup_{i=1}^t\Ass_{R_\p}\Ext_{R_\p}^n((M_i^\dag)_\p/I(M_i^\dag)_\p,\omega_{R_{\p}})$.
Hence the ideal $\p R_\p$ is in $\bigcup_{i=1}^t\Ass_{R_\p}\Ext_{R_\p}^n((M_i^\dag)_\p/I(M_i^\dag)_\p,\omega_{R_{\p}})$, and therefore $\p$ is in $\bigcup_{i=1}^t\Ass_R\Ext_R^n(M_i^\dag/IM_i^\dag,\omega_R)$.
This completes the proof of the theorem.
\end{proof}

\section{Appendix: Discreteness of F-jumping exponents}
Here, using arguments similar to those in Section 2, we study the discreteness of F-jumping exponents of ideals of rings with finite F-representation type. 

Let $R$ be a Noetherian reduced ring of characteristic $p > 0$ and $M$ an $R$-module. 
Denote by $R^{\circ}$ the set of elements of $R$ which are not in any minimal prime ideal. 
For each integer $e \ge 0$, we denote $\F^e(M) = \F_R^e(M) := M \otimes_R {}^e\! R$ and regard it as an $R$-module by the action of $R$ on ${}^e\! R$ from the right. 
Then we have the induced $e$-times iterated Frobenius map $F_M^e \colon M \to \F^e(M)$. 
The image of $z \in M$ via this map is denoted by $z^q:= F_M^e(z) 
\in \F^e(M)$, where $q=p^e$. 

Now we recall the definition of $\a^t$-tight closure. 
\begin{defn}[\textup{\cite[Definition 6.1]{HY}}]
Let $\a$ be an ideal of a Noetherian reduced ring $R$ of characteristic $p>0$ such that $\a \cap R^{\circ} \ne \emptyset$ and let $t>0$ be a real number. 
Let $M$ be a (not necessarily finitely generated) $R$-module. 
The $\a^t$-tight closure of the zero submodule in $M$, denoted by $0_M^{*\a^t}$, is defined to be the submodule of $M$ consisting of all elements $z \in M$ for which there exists $c \in R^{\circ}$ such that 
$$c \a^{\lceil tq \rceil} z^q=0 \in \F^e(M)$$
for all large $q=p^e$. 
\end{defn} 

\begin{defn}[\textup{cf.\! \cite[Definition 1.6, Proposition-Definition 1.9]{HY}}]\label{taudef}
Let $\a$ be an ideal of an excellent Noetherian reduced ring $R$ of characteristic $p>0$ such that $\a \cap R^{\circ} \ne \emptyset$ and let $t >0$ be a real number. 
Let $E =\bigoplus_{\m} E_R(R/\m)$ be the direct sum, taken over all maximal ideals $\m$ of $R$, of the injective hulls of the residue fields $R/\m$. 
\begin{enumerate}
\renewcommand{\labelenumi}{\textup{(\roman{enumi})}}
\item
The \textit{generalized test ideal} ${\tau}(\a^t)$ of $\a$ with exponent $t$ is defined to be 
$${\tau}(\a^t)=\tau(R,\a^t)=\displaystyle\bigcap_{M\subset E} \Ann_R(0^{*\a^t}_M),$$
where $M$ runs through all finitely generated $R$-submodules of $E$. 

\item The ideal $\widetilde{\tau}(\a^t)$ is defined to be
$$\widetilde{\tau}(\a^t)=\widetilde{\tau}(R,\a^t)=\Ann_R(0^{*\a^t}_E).$$

\item We say that an element $c \in R^{\circ}$ is an \textit{$\a^t$-test element} for $E$ if $c\a^{\lceil tq \rceil}z^q=0$ in $\F^e(E)$ for all $q=p^e$ whenever $z \in 0^{*\a^t}_E$. 
\end{enumerate}
\end{defn}

\begin{prop}[\textup{cf.\! \cite[Theorem 1.7, Theorem 1.13]{HY}, \cite[Theorem 3.3]{LS}}]\label{test result}
Let $R$ be an excellent Noetherian reduced ring of characteristic $p>0$ and $\a \subseteq R$ be an ideal such that $\a \cap R^{\circ} \ne \emptyset$.  
\renewcommand{\labelenumi}{\textup{(\arabic{enumi})}}
\begin{enumerate}
\item Assume that one of the following conditions holds: \textup{(i)} $R$ is a $\Q$-Gorenstein normal local ring; \textup{(ii)} $R=\bigoplus_{n \ge 0}R_n$ is  a graded ring with $R_0$ a field and $\a$ is a homogeneous ideal.
Then for all real numbers $t >0$, 
$${\tau}(\a^t)=\widetilde{\tau}(\a^t).$$  

\item Let $c \in R^{\circ}$. 
If $R$ is F-finite and the localized ring $R_c$ is strongly F-regular, then some power $c^n$ of $c$ is an $\a^t$-test element for $E$ for all ideals $\a \subseteq R$ such that $\a \cap R^{\circ} \ne \emptyset$ and for all real numbers $t>0$. 
\end{enumerate}
\end{prop}

Let $R$ be an F-finite Noetherian reduced ring of characteristic $p>0$. We say that $R$ is \textit{strongly F-regular} if for every $c \in R^{\circ}$, there exists $q=p^e$ such that $c^{1/q}R \hookrightarrow R^{1/q}$ splits as an  $R$-module homomorphism.  When the ring is strongly F-regular, we can describe the generalized test ideal ${\tau}(\a^t)$ in a simple form. 
\begin{prop}\label{test ideal}
Let $R$ be a strongly F-regular ring of characteristic $p>0$ and $\a \subseteq R$ be an ideal such that $\a \cap R^{\circ} \ne \emptyset$.  
\begin{enumerate}
\renewcommand{\labelenumi}{\textup{(\arabic{enumi})}}
\item
$$\widetilde{\tau}(\a^t)=\bigcup_{q=p^e} I_e(\a^{\lceil tq \rceil}, R)$$
$($see Definition \ref{idealdef} for the definition of the ideal $I_e(\a^{\lceil tq \rceil}, R)$ of $R)$.  
Since $\{I_e(\a^{\lceil tq \rceil},R)\}_{q=p^e}$ forms an ascending chain of ideals in $R$, this union stabilizes after finitely many steps. 
In particular, the ideal $\widetilde{\tau}(\a^t)$ is equal to $ I_e(\a^{\lceil tq \rceil}, R)$ for all sufficiently large $q=p^e$. 
\item
If $R$ is $\Q$-Gorenstein, then for all sufficiently large $q=p^e$, 
$$\tau(\a^t)=\widetilde{\tau}(\a^t)=I_e(\a^{\lceil tq \rceil},R).$$
\end{enumerate}
\end{prop}

\begin{proof}
First, we will prove that $\widetilde{\tau}(\a^t) \supseteq I_e(\a^{\lceil tp^e \rceil}, R)$ for all integers $e \ge 0$. 
By Proposition \ref{test result} (2), the unit element $1$ is an $\a^t$-test element for $E$. 
It means that $\a^{\lceil tp^e \rceil}z^{p^e}=0$ in $E \otimes_R {}^e\!R$ for all $z \in 0^{*\a^t}_E$. 
An $R$-linear map $\varphi: {}^e\!R \to R$ induces an $R$-linear map $E \otimes_R {}^e\!R \to E$ sending $\a^{\lceil tp^e \rceil}z^{p^e}$ to $\varphi(\a^{\lceil tp^e \rceil})z$. 
Thus, $\varphi(\a^{\lceil tp^e \rceil})0^{*\a^t}_E=0$ in $E$, which implies that $\varphi(\a^{\lceil tp^e \rceil}) \subseteq \widetilde{\tau}(\a^t)$. 

To show the reverse inclusion, it suffices to prove it after passing to localization at a maximal ideal $\m$ in $R$. 
Since  the formation of the ideal $I_e(\a^{\lceil tp^e \rceil}, R)$ obviously commutes with localization,  
$$\left(\bigcup_{e \ge 0} I_e(\a^{\lceil tp^e \rceil}, R)\right)R_{\m} = \bigcup_{e \ge 0} I_e(\a_{\m}^{\lceil tp^e \rceil}, R_{\m}).$$
Also, by \cite[Lemma 2.1]{HT}, we have 
$$\widetilde{\tau}(R_{\m}, \a_{\m}^t)=\bigcup_{e \ge 0} I_e(\a_{\m}^{\lceil tp^e \rceil}, R_{\m}).$$
Thus, combining the above two equalities, we obtain an inclusion
$$\widetilde{\tau}(R,\a^t)R_{\m} \subseteq \widetilde{\tau}(R_m,\a_m^t)=\bigcup_{e \ge 0} I_e(\a_{\m}^{\lceil tp^e \rceil}, R_{\m})=\left(\bigcup_{e \ge 0} I_e(\a^{\lceil tp^e \rceil}, R)\right)R_{\m}.$$

The latter assertion follows from a similar argument, because the generalized test ideal $\tau(R_{\m},\a_{\m}^t)$ coincides with the ideal $\widetilde{\tau}(R_{\m}, \a_{\m}^t)$ when $R$ is normal $\Q$-Gorenstein. 
\end{proof}

The family of ideals $\{\widetilde{\tau}(\a^t)\}_{t \ge 0}$ is right continuous in $t$. 
\begin{lem}[\textup{cf. \!\cite[Corollary 2.16]{BMS}}]\label{conti}
Let $R$ be a strongly F-regular ring of characteristic $p>0$. 
For every ideal $\a \subseteq R$ such that $\a \cap R^{\circ} \ne \emptyset$ and for every real number $t > 0$, there exists $\epsilon >0$ such that $\widetilde{\tau}(\a^t)=\widetilde{\tau}(\a^{t'})$ for all $t' \in [t, t+\epsilon)$. 
\end{lem}

\begin{proof}
First we will prove the following claim. 
\begin{cl}
For every $e \ge 0$, there exists $e' \ge e$ such that $I_e(\a^{\lceil tp^e \rceil}, R) \subseteq I_{e'}(\a^{\lceil tp^{e'} \rceil+1}, R)$. 
\end{cl}
\begin{proof}[Proof of Claim]
Let $c \in \a \cap R^{\circ}$. 
Since $R$ is strongly F-regular, there exist $q_0=p^{e_0}$ and an $R$-linear map $\psi:{}^{e_0}\!R \to R$ sending $c$ to $1$. 
For any $\varphi_e \in \Hom_R({}^e\!R, R)$, $\psi$ induces an $R$-linear map $\varphi_{e+e_0}:{}^{e+e_0}\!R \to R$ sending $c$ to $\varphi_e(1)$ so that 
$$\varphi_e(\a^{\lceil t q \rceil})=\varphi_{e+e_0}(c(\a^{\lceil t q \rceil })^{[q_0]}) \subseteq I_{e+e_0}(\a^{\lceil tqq_0 \rceil+1}, R).$$ 
\end{proof}

Take sufficiently large $e$ so that $\widetilde{\tau}(\a^t)=I_e(\a^{\lceil tp^e \rceil}, R)$ by Proposition \ref{test ideal}.  Then by the above claim, there exists $e' \ge e$ such that $I_e(\a^{\lceil tp^e \rceil}, R) \subseteq I_{e'}(\a^{\lceil tp^{e'} \rceil+1}, R)$. 
Put $t'=t+1/p^{e'}$. Then, by Proposition \ref{test ideal} again,  the above inclusion implies  the required inclusion
$$\widetilde{\tau}(\a^t) \subseteq I_{e'}(\a^{\lceil t'p^{e'} \rceil}, R) \subseteq \widetilde{\tau}(\a^{t'}).$$
This completes the proof. 
\end{proof}

Thanks to Proposition \ref{test result} (1) and Lemma \ref{conti}, one can define jumping exponents for the generalized test ideals $\tau(\a^t)$.   
\begin{defn}[\textup{\cite[Definition 2.17]{BMS}}]
Let $\a$ be an ideal of an excellent Noetherian reduced ring $R$ of characteristic $p>0$ such that $\a \cap R^{\circ} \ne \emptyset$. 
Assume that one of the following conditions holds: \textup{(i)} $R$ is a $\Q$-Gorenstein normal local ring; \textup{(ii)} $R=\bigoplus_{n \ge 0}R_n$ is  a graded ring with $R_0$ a field and $\a$ is a homogeneous ideal.
Then a positive real number $t$ is said to be an \textit{F-jumping exponent} of $\a$ if $\tau(\a^t) \subsetneq \tau(\a^{t-\epsilon})$ for every $\epsilon>0$. 
\end{defn}

\begin{lem}[\textup{cf. \!\cite[Proposition 3.2]{BMS}}]\label{degree}
Let $R=\bigoplus_{n \ge 0}R_n$ be a strongly F-regular graded ring with $R_0$ a field of characteristic $p>0$ and let $\a \subseteq R$ be a homogeneous ideal such that $\a \cap R^{\circ} \ne \emptyset$. 
Assume in addition that $R$ has finite graded F-representation type. 
Then there exists an integer $r \ge 0$ depending only on $R$ such that 
if $\a$ is generated by elements of degree at most $d$, then for every $t>0$, 
the generalized test ideal ${\tau}(\a^t)$ can be generated by elements of degree at most $td+r$. 
\end{lem}
\begin{proof}
Since $R$ has finite graded F-representation type, there exist finitely generated $\Q$-graded $R$-modules $M_1, \dots, M_s$ such that for each $e \in \N$, 
$${}^e\!R \cong \bigoplus_{i=1}^s \bigoplus_{j=1}^{n_{ei}}M_i(\alpha_{ij}^{(e)})$$
as a $\Q$-graded R-module. 
Without loss of generality, we may assume that the grading on the $M_i$ is normalized so that $[M_i]_0 \ne 0$ and $[M_i]_{\alpha}=0$ for any $\alpha<0$. 
Then by the claim in the proof of Theorem \ref{stabilize}, there exists a nonnegative rational number $C_1$ such that $-\alpha_{ij}^{(e)} \le C_1$ for all $e \in \N$, $i=1, \dots, s$ and $j=1, \dots, n_{ei}$.  

Now fix any $e \in \N$. The ideal $\a^{\lceil tp^e \rceil}$ is generated by elements of degree at most $\lceil tp^e \rceil d$. Choose such generators $g_1, \dots, g_{a_{e}}$ for $\a^{\lceil tp^e \rceil}$. 
For each $k=1, \dots, a_e$, write $g_k=\sum_{i=1}^s\sum_{j=1}^{n_{ei}}m_{ijk}^{(e)}$ with $m_{ijk}^{(e)} \in M_i(\alpha_{ij}^{(e)})$ via the identification of ${}^e\!R$ with $\bigoplus_{i=1}^s \bigoplus_{j=1}^{n_{ei}}M_i(\alpha_{ij}^{(e)})$. Then one has 
$$\deg_{M_i(\alpha_{ij}^{(e)})}(m_{ijk}^{(e)}) \le \deg_{{}^e\!R}(g_k)=\frac{1}{p^e}\deg_{R}(g_k)  \le \frac{\lceil tp^e \rceil d}{p^e}.$$
Here we take a rational number $C_2$ such that  for every $i=1, \dots, s$, $\Hom_{R}(M_i, R)$ can be generated by elements of degree at most $C_2$. 
Then 
\begin{align*}
I_e(\a^{\lceil tp^e \rceil}, R)
&=\sum_{i=1}^s\sum_{j=1}^{n_{ei}}\sum_{k=1}^{a_e}\Hom_{R}(M_i(\alpha_{ij}^{(e)}), R) \cdot m_{ijk}^{(e)}\\
&\subseteq \sum_{i=1}^s\sum_{j=1}^{n_{ei}}\Hom_{R}(M_i(\alpha_{ij}^{(e)}), R) \cdot [M_i(\alpha_{ij}^{(e)})]_{\le \frac{\lceil tp^e \rceil d}{p^e}}\\
&\subseteq \sum_{i=1}^s \Hom_{R}(M_i, R) \cdot [M_i]_{\le \frac{\lceil tp^e \rceil d}{p^e}+C_1}\\
& \subseteq [R]_{\le  \frac{\lceil tp^e \rceil d}{p^e}+C_1+C_2}.
\end{align*}
Thus, given $t$, if $e$ is enough large, then 
$I_e(\a^{\lceil tp^e \rceil}, R)$ can be generated by elements of degree at most  $\lfloor td \rfloor +C_1+C_2$. 
On the other hand, by Propositions \ref{test result} (1) and \ref{test ideal} (1), $I_e(\a^{\lceil tp^e \rceil}, R)$ coincides with the ideal $\tau(\a^t)$ for sufficiently large $e$. 
\end{proof}

Using Lemma \ref{degree}, we can generalize the result due to Blickle, Musta\c t\v a, Smith (\cite[Theorem 3.1 (1)]{BMS}) to the case of rings with finite F-representation type. 
\begin{thm}\label{discrete}
Let $R=\bigoplus_{n \in \N}R_n$ be a strongly F-regular graded ring with $R_0$ a field of characteristic $p>0$ and let $\a$ be a homogeneous ideal of $R$ such that $\a \cap R^{\circ} \ne \emptyset$. 
Assume in addition that $R$ has finite graded F-representation type. 
Then the set of F-jumping exponents of $\a$ have no accumulation points. 
\end{thm}

\begin{proof}
Once we accept Lemmas \ref{conti} and \ref{degree},  the proof is quite similar to that of \cite[Theorem 3.1]{BMS}. 

Suppose that we have a sequence of F-jumping exponents $\{\alpha_m\}_{m}$ of $\a$ having a finite accumulation point $\alpha$. By virtue of Lemma \ref{conti}, we can see that $\alpha_m <\alpha$ for sufficiently large $m$. After replacing this sequence by a subsequence, we may assume that $\alpha_m <\alpha_{m+1}$ for every $m$. 
Assume now that $\a$ is generated by elements of degree at most $d$. Then, by Lemma \ref{degree}, there exists an integer $r$ (depending only on $R$) such that every generalized test ideal $\tau(\a^{\alpha_m})$ is generated by elements of degree at most $\alpha d +r$. 
Since $[R]_{\le \alpha d+r}$ is a finite-dimensional vector space over $R_0$, 
 the descending chain of subspaces in $[R]_{\le \alpha d+r}$
$$\tau(\a^{\alpha_1})\cap [R]_{\le \alpha d+r} \supseteq \tau(\a^{\alpha_2})\cap [R]_{\le \alpha d+r} \supseteq \tau(\a^{\alpha_3}) \cap [R]_{\le \alpha d+r} \supseteq \dots$$
stabilizes. Since $\tau(\a^{\alpha_m})\cap [R]_{\le \alpha d+r}$ generates $\tau(\a^{\alpha_m})$ as an $R$-module for every $m$, this implies that the descending chain of generalized test ideals 
$$\tau(\a^{\alpha_1}) \supsetneq \tau(\a^{\alpha_2}) \supsetneq \tau(\a^{\alpha_3}) \supsetneq \dots$$ stabilizes. This gives a contradiction. 
\end{proof}

The rationality of F-jumping exponents is a little more subtle. 
\begin{ques}
Let the notation be the same as in Theorem \ref{discrete}. 
Does there exist a power $q_0=p^{e_0}$ of $p$ such that if $\alpha$ is an F-jumping exponent of $\a$, then so is $q_0 \alpha$?
\end{ques}

If the answer to the above question is yes, then we can prove the rationality of F-jumping exponents of $\a$ under the same assumption as that of Theorem \ref{discrete}. 

\begin{rem}
Let $R=\bigoplus_{n \ge 0}R_n$ be a strongly F-regular graded ring with $R_0$ a field of characteristic $p>0$. 
If $R$ is a $\Q$-Gorenstein ring and the order of the canonical class in the divisor class group $\mathrm{Cl}(R)$ is not divisible by $p$, then the set of F-jumping exponents is discrete and all F-jumping exponents are rational numbers (without assuming that $R$ has finite graded F-representation type). 
We refer to \cite{HuT} for details. 
\end{rem}

\begin{small}
\begin{acknowledgement}
The authors are indebted to Mircea Musta{\c{t}}{\v{a}} for pointing out a mistake in a previous version of this paper and to Manuel Blickle for suggesting to include a result on F-jumping exponents. 
They are also grateful to Holger Brenner, Shiro Goto, Mel Hochster, Yuji Kamoi, Ken-ichiro Kawasaki,  Kazuhiko Kurano, Takafumi Shibuta, Anurag Singh, Karen Smith, Kei-ichi Watanabe, Yongwei Yao and Ken-ichi Yoshida for valuable conversations and helpful suggestions. 
They were partially supported by Grant-in-Aid for Young Scientists (B) 17740021 and 19740008 from JSPS, respectively. 
The first author was also partially supported by Program for Improvement of Research Environment for Young Researchers from SCF commissioned by MEXT of Japan.
\end{acknowledgement}
\end{small}

\end{document}